\documentclass[11pt]{amsart}  
\usepackage{amssymb,amsmath,amsthm,amscd}
\usepackage{graphicx}
\numberwithin{equation}{section}

\newtheoremstyle{fancy1}{10pt}{10pt}{\itshape}{12pt}{\textsc\bgroup}{.\egroup}{8pt}{
}
\newtheoremstyle{fancy2}{10pt}{10pt}{}{12pt}{\itshape}{.}{8pt}{ }

\theoremstyle{fancy1}

\newtheorem{lem}[equation]{Lemma}

\newtheorem{thm}[equation]{Theorem}
\newtheorem{problem}{Problem}

\newtheorem*{main*}{Theorem}
\newtheorem*{conj}{Conjecture}
\newtheorem*{cor*}{Corollary}
\newtheorem*{prop*}{Proposition}

\newtheorem*{problem*}{Problem}

\theoremstyle{fancy2}

\newtheorem*{rem*}{Remark}

\newcommand{\cref}[1]{Corollary~\ref{#1}}

\newcommand{\RP}{\mathbb{R\mkern1mu P}}
\newcommand{\CP}{\mathbb{C\mkern1mu P}}
\newcommand{\HP}{\mathbb{H\mkern1mu P}}
\newcommand{\CaP}{\mathrm{Ca}\mathbb{\mkern1mu P}^2}
\newcommand{\Sph}{\mathbb{S}}
\newcommand{\Disc}{\mathbb{D}}

\newcommand{\Z}{{\mathbb{Z}}}
 
\newcommand{\N}{{\mathbb{N}}}

\newcommand{\G}{\ensuremath{\operatorname{G}}}

\newcommand{\SO}{\ensuremath{\operatorname{SO}}}

\newcommand{\Sp}{\ensuremath{\operatorname{Sp}}}
\newcommand{\U}{\ensuremath{\operatorname{U}}}
\newcommand{\SU}{\ensuremath{\operatorname{SU}}}

\newcommand{\T}{\ensuremath{\operatorname{T}}}
\renewcommand{\S}{\ensuremath{\operatorname{S}}}

\def\con#1=#2(#3){#1 \equiv #2 \bmod{#3}}

\DeclareMathOperator{\symrank}{symrank}
\DeclareMathOperator{\symdeg}{symdeg}
\DeclareMathOperator{\cohom}{cohom}

\newcommand{\no}{\noindent}

\begin{document}

\title{Developments around positive sectional curvature}

\author{Karsten Grove}
\address{University of Notre Dame}
\email{kgrove2@nd.edu}

\thanks{ Supported in part by  a grant from the National Science
Foundation. }

\maketitle

\begin{abstract}
This is not in any way meant to be a complete survey on positive curvature. Rather it is a short essay on the fascinating changes  in the landscape surrounding positive curvature. In particular, details and many results and references are not included, and things are not presented in chronological order.
\end{abstract}

\bigskip

Spaces of positive curvature have always enjoyed a particular role in Riemannian geometry. Classically, this class of spaces form a natural and vast extension of \emph{spherical geometry}, and in the last few decades their importance for the study  of general manifolds with a lower curvature bound via \emph{Alexandrov geometry} has become apparent.

 The importance of Alexandrov geometry to Riemannian geometry stems from the fact that there are several natural geometric operations that are closed in Alexandrov geometry but not in Riemanian geometry. These include \emph{taking Gromov Hausdorff limits}, \emph{taking quotients} by isometric group actions, and \emph{forming joins} of positively curved spaces. In particular, limits (or quotients) of Riemanian manifolds with a lower (sectional) curvature bound are Alexandrov spaces, and only rarely Riemannian manifolds. Analyzing limits frequently involves blow ups leading to spaces with \emph{non-negative curvature} as, e.g., in Perelman's work on the geometrization conjecture. Also the infinitesimal structure of an Alexandrov space is expressed via its ``tangent spaces", which are cones on positively curved spaces. Hence the collection of all compact positively curved spaces (up to scaling) agrees with the class of all possible so-called \emph{spaces of directions},  in Alexandrov spaces. So spaces of positive curvature play the same role in Alexandrov geometry as round spheres do in Riemannian geometry.  
 
 In addition to positively, and nonnegativey curved spaces, yet another class of spaces has emerged in the general context of convergence  under a lower curvature bound, namely \emph{almost nonnegatively curved} spaces. These are spaces allowing metrics with diameter say $1$, and lower curvature bound arbitrarily close to $0$. They are expected to play a role among spaces with a lower curvature bound, analogous to that almost flat spaces play for spaces with bounded curvature. In summary, the following classes of spaces play essential roles in the study of spaces with a lower curvature bound:
 
 \begin{center}
 $ \mathcal{P} \subset  \mathcal{P}_0 \subset \mathcal{P}_{0^{-}}$
 \end{center} 
 
 \no corresponding to \emph{positively curved}, \emph{nonnegatively curved}, and \emph{almost nonnegatively curved} spaces. Here the first and the last class enjoy the useful property of being stable under small perturbations. Among all manifolds, these form ``the tip of the iceberg". Yet, aside from being  \emph{Nilpotent spaces} (up to finite covers) and having a priory \emph{bounded topology} in terms of generators for homology- and fundamental groups, only a few general obstructions are known, and none in the simply connected case. Moreover, so far only obstructions on fundamental groups distinguish the three classes.
 \bigskip
 
 The study of the two first classes has also played a significant role in the development of \emph{comparison theorems}, and in this way also influenced the general development of Riemannian geometry. Much of this work even originated in connection with the desire to characterize spheres among spaces with positive curvature, resulting in so-called \emph{sphere theorems}. It is indeed remarkable how many of the known general tools  have been developed during proofs of sphere theorems. This not only applies to a number of so-called \emph{comparison theorems}, like the \emph{Rauch Comparison Theorem} for \emph{Jacobi fields} and the \emph{Toponogov Comparison Theorem} for \emph{geodesic triangles}, but also includes \emph{critical point theory for distance functions} as well as the \emph{Ricci flow}.
 
   \bigskip
 
 Throughout the rest of this article our main focus will be on manifolds with positive (sectional) curvature. Except for the beautifully rich and simple trick provided by Synge,  there have been three general approaches used in an attempt to gain understanding of the class, all starting in some sense with the sphere as a uniquely determined extremal object. In short these approaches have been guided by \emph{Shape}, \emph{Size}  and \emph{Symmetry}. Here the first two of these have played an important role in the development of the tools alluded to above and hence also to the known obstructions. The latter, on the other hand has provided a natural framework for the discovery and construction of examples, an area of pivotal importance for the subject.  We will discuss each of them below.
 
 \smallskip
 
 The reader should also consult two recent survey articles on manifolds with positive and nonnegative curvature by Wilking \cite{Wi} and Ziller \cite{Z1}. Most of what we will discuss here can be found in one of these surveys where more details are given as well as references to original papers. Other surveys that cover topics mentioned in this essay are Petersen \cite{Pe} and Rong \cite{Ro}. The surveys by Plaut \cite{Pl} and Petrunin \cite{Pt} provide good access to Alexandrov spaces. Other related surveys are \cite{Z2} and \cite{G1,G2}. We will use a non-traditional approach and only provide other references if they cannot easily be found in one of these sources, or when the reference therein was incomplete.
 
\section{Structure}

Without stating it again, we assume throughout that all manifolds considered are \emph{complete} Riemannian manifolds.

\bigskip

The following fundamental result due to Cheeger and Gromoll provides an amazingly simple structure of all non-compact manifolds in $\mathcal{P}_0$:

\begin{thm}[Soul Theorem]
Any non-compact manifold $M \in \mathcal{P}_0$  is diffeomorphic to the total space $E$ of a vector bundle over a compact  manifold $S \in \mathcal{P}_0$. In fact, any  soul $S$ is a totally geodesic submanifold of $M$ and $E$ is the normal bundle to $S$ in $M$.
\end{thm}

\no Except for the fundamental:

\begin{problem}
Which vector bundles over a compact $S \in \mathcal{P}_0$ admit a metric with nonnegative curvature?
\end{problem}

\no the soul theorem to a large extent reduces the topological study of nonnegatively curved manifolds to compact ones. So far, obstructions in this problem are known only for certain manifolds $S$ with infinite fundamental group by work of Ozaydin-Walschap and  of Belegradek-Kapovitch. The problem is difficult even when $S$ is a sphere. Here work of Grove-Ziller shows that all bundles have such metrics when the dimension of the base is at most 5.

For positive curvature, one has a complete answer due to Gromoll and Meyer (known prior to the soul theorem):

\begin{thm}
Any non-compact manifold $M \in \mathcal{P}$  is diffeomorphic to euclidean space.
\end{thm}

\no In the language of the soul theorem, a soul in this case is simply a point. As shown by Perelman this remains true in the more general case where $M$ has positive curvature on an open set. The key geometric tools in all of this work are the Toponogov comparison theorem and a deep study of convex sets and concavity properties of distance functions in nonnegative curvature.

\bigskip

In the remaining part of this essay, we will confine our discussion to compact manifolds.

\bigskip

The first global theorems about positive curvature were special applications of the \emph{Second Variation Formula for Arc Length}, and hence in essence \emph{Index Theorems} in the context of \emph{Morse Theory} for suitable path spaces. The particular method pioneered by Synge (variations in directions of parallel fields) applies in general situations where \emph{boundary conditions} for the relevant path spaces are \emph{totally geodesic}. The original application deals with \emph{orientation}:

\begin{thm}[Orientation]
Any odd dimensional  manifold in $\mathcal{P}$ is orientable, and an orientable even dimensional manifold in $\mathcal{P}$ is simply connected.
\end{thm}

A recent application of the Synge trick due to Wilking illustrates how severely the presence of totally geodesic submanifolds (with small codimension) restricts the structure of manifolds in $\mathcal{P}$:

\begin{lem}[Connectivity]
The inclusion map of a totally geodesic submanifold $V^{n-k} \subset M^n \in \mathcal{P}$ is $n-2k+1$ connected.
\end{lem}

Recall here that a map $f: X \to Y$ is $\ell$-connected if the induced map $f_q: \pi_q(X) \to \pi_q(Y)$ on homotopy groups is onto for $q=\ell$ and an isomorphism for $q < \ell$. 

For $k=1$ a convexity argument as in the soul theorem immediately implies that such a manifold $M$ is homeomorphic to a sphere or is double covered by a manifold homeomorphic to a sphere. The following is a natural open

\begin{problem}
Is any $M \in \mathcal{P}$ with a codimension one totally geodesic submanifold diffeomorphic to either the standard sphere or the real projective space?
\end{problem}

Also for $k=2$, the following could have a positive answer:

\begin{problem}
Is any (simply connected) $M \in \mathcal{P}$ with a codimension two totally geodesic submanifold diffeomorphic to either the standard sphere or the complex projective space?
\end{problem}

From now on we will use the terminology CROSS for a simply connected  \emph{compact rank one symmetric space}, i.e., such a manifold is either the sphere $\Sph^n$, the complex-, or quaternionic- projective space $\CP^n$, or $\HP^n$, or the Cayley plane $\CaP$ with their canonical metrics. Each type forms a \emph{chain} (with fixed codimension), i.e., an infinite sequence of manifolds  $M_i$, with $M_i \subset M_{i+1}$ isometrically and in particularly  totally geodesically. Positively curved chains are clearly severely restricted by the connectivity lemma. In fact

\begin{problem}
Are the CROSS's the only positively curved chains?
\end{problem}

It should be pointed out that at present simply connected manifolds $M \in \mathcal{P}$ not diffeomorphic to a CROSS are known only in dimensions 6 (two), 7 (infinitely many), 12 (one), 13 (infinitely many), and 24 (one).

\section{Symmetry}

The round sphere, $\Sph^n$ is  characterized as the simply connected closed manifold with maximal symmetry degree, i.e., with the largest dimensional Lie group acting on it. This naturally leads to questions of describing manifolds with large transformation groups even without any curvature assumptions as, e.g., in work of W.-Y. Hsiang and many others. It also suggest the following approach to investigate manifolds of positive curvature:

\smallskip

\begin{center}
\emph{Desribe the structure and ultimately classify positively
 curved manifolds\\ with large isometry groups}
\end{center}

\smallskip

One of the motivations for this, is that it provides a systematic guide towards the discovery and construction of new examples, a task which is known to be very difficult. In retrospect, the classification of (simply connected) homogeneous manifolds with positive curvature, due to the combined efforts of Berger, Wallach, Aloff-Wallach and Berard-Bergery, is an example of this. Non-CROSS'es occur in dimensions 6 ($W^6$), 7 ($B^7$ and $A^7_{p,q}, p,q \in \Z$), 12 ($W^{12}$), 13 ($B^{13}$), and 24 ($W^{24}$). The driving force in that work, however, was based on the fact that Lie groups with biinvariant metrics have nonnegative curvature and that taking quotients only increases curvature. The same method was behind the construction due to Eschenburg in dimensions 6 ($E^6$) and 7 ($E^7_{k,l}, k, l \in \Z^3$) and later Bazaikin in dimension 13 ($B^{13}_{q}, q \in \Z^5$) of so-called \emph{biquotients} with positive curvature. In both cases curvature computations are essentially reduced to Lie algebra computations. 

\smallskip

The same type of problems are natural for the other classes $ \mathcal{P}_0 \subset \mathcal{P}_{0^{-}}$ as well, and rich classes in  $\mathcal{P}_0$, including many exotic spheres have been discovered in this way.

\bigskip

There are several natural \emph{measures for largeness}. The most obvious ones are: \emph{symdeg, symrank, cohom}, i.e., the degree of symmetry = dimension of Isometry group, the symmetry rank = rank of isometry group, the cohomogeneity = dimension of orbit space. Another very useful measurement, \emph{cofix} takes into account the fixed point set and essentially measures cohomogeneity modulo fixed point sets, called the \emph{fixed point cohomogeneity} = cohom - dim fix - 1. For each of these invariants there are partial or complete solutions to the program described above. 

\bigskip

In all cases, orbit spaces play a significant role. These orbit spaces are Alexandrov spaces with positive curvature, and this vaguely prevents ``too many", ``too singular" points. Typically, the orbit space will have non-empty boundary, and ``soul theorems'' for these spaces become powerful tools. On the manifold itself, the existence of \emph{totally geodesic submanifolds} as fixed point sets of groups of isometries are common in this context, often leading to restrictions via the connectivity lemma of section one.

\bigskip

Torus actions play a particular role for several reasons. On the one hand, it follows from the Cheeger-Fukaya-Gromov theory, that collapse with bounded curvature of simply connected manifolds is given in terms of such actions. Secondly, for manifold in $\mathcal{P}$ the method of Synge implies that an isometric action by a torus $\T$ has fixed points in even dimensions, and either fixed points or circle orbits in odd dimensions. In particular, there are isotropy groups of maximal rank in even dimension and at most corank 1 in odd dimensions.

\bigskip

 For the symmetry rank one has the following rigidity and ``pinching" theorems due to Grove-Searle and Wilking respectively:
 
 \begin{thm}[Rank Rigidity]
Any $M^n \in  \mathcal{P}$ spectacularly has $\symrank (M) \le [\frac{n+1}{2}]$ with equality if and only if $M$ is diffeomorphic to either $\Sph^n$, $\CP^{n/2}$, or a lense space $\Sph^n/ \Z_k$.
\end{thm}

\no Moreover, the actions are known as well.

\begin{thm}[Rank Pinching]
Any positively curved simply connected $n$-manifold $M$ with $\symrank (M) \ge n/4 + 1$ and $n \ne 7$ is homotopy equivalent to a CROSS.
\end{thm}

\no In this formulation, results have been combined: It is due to Wilking for $n \ge 10$ , to Fang-Rong for $7 < n < 10$ and is covered by the rigidity theorem in dimensions  $n <  7$ .  The  conclusion fails spectacularly in dimension 7. In fact all known examples, i.e., all positively curved  Eschenburg and Aloff Wallach spaces have isometry group of rank 3.

When  $\symrank (M) =  [\frac{n+1}{2}] - 1$, we say that $M$ has \emph{almost maximal} symmetry rank. The only dimensions not covered by the above result are dimensions 4, 5, 6, and 7 where this is 1, 2, 2, and 3. Here, the flag manifold $W^6 = \SU(3)/\T^2$ and the Eschenburg ``flag" $E^6 = \SU(3)//\T^2$ both have almost maximal symmetry rank. Since dimensions 4 and 5 have been resolved (see below), this suggests the following 

\begin{problem}
Classify $M^n \in \mathcal{P}$ with almost maximal symmetry rank when $n = 6$ and 7.
\end{problem}

\no As indicated, all known examples in $ \mathcal{P}$ in these dimensions have almost maximal symmetry rank.

In lower dimensions there is a complete classification:

\begin{thm}[Strong Rank Pinching]
Any  $M^n \in \mathcal{P}$ with almost maximal symmetry rank and $n \le 5$ is diffeomorphic to a  CROSS.
\end{thm}

In dimensions 2 and 3, of course, no symmetry assumptions are  needed. In dimension 5, the result is due to Rong and since there are no exotic spheres in this dimension it suffices to prove that such a manifold is a homotopy sphere. The \emph{topological classification} in dimension 4 is due to Hsiang and Kleiner. By Freedman's topological classification of simply connected 4-manifolds, it suffices to prove that the Euler characteristic of such a manifold is at most 3. Since this is also the Euler characteristic of the fixed point set of the circle action, the result follows from a geometric analysis of the fixed point set. 

When confined to smooth circle actions (cf. \cite{FSS}), the extensive work of Fintushel \cite{Fi} and Pao \cite{Pa} gives, modulo the Poincar\'{e} conjecture, in particular a complete classification of such actions on simply connected 4-manifolds, in terms of so-called \emph{weighted} orbit space data. Thus, the combination of the work of Hsiang-Kleiner, Fintushel and Perelman provides a proof of the above theorem in dimension 4. Extending the fixed point analysis to non-negative curvature as in Kleiner's thesis and in the work of Searle-Yang, and again combining this with the work of Fintushel and Perelman immediately yields the following classification:

\begin{thm}
Any non-negatively curved simply connected 4-manifold with an isometric circle action is diffeomorphic to one of
 $\Sph^4, \CP^2, \Sph^2 \times \Sph^2$,
 $ \CP^2 \# \pm \CP^2$.
\end{thm}

\no A classification up to equivariant diffeomorphism is not yet known, since a full understanding of the possible weighted orbit space data in nonnegative and positive curvature has not yet been achieved. In the case of positive curvature, however, we propose the following

\begin{conj}
An isometric circle action on a simply connected positively curved 4-manifold is equivariently diffeomorphic to a linear circle action on either  $\Sph^4$  or $ \CP^2$.
\end{conj}

The theorem above of course provides support for the classical Hopf conjecture, that there is no metric on $\Sph^2 \times \Sph^2$ with positive curvature. Indeed if there is one, it can have at most a finite isometry group. One can speculate that the list in the theorem are all simply connected manifolds $\mathcal{P}_0$. Topologically this would follow from the ellipticity conjecture (see ``Size" section).

\smallskip

Replacing the torus $\T^k$ by any compact Lie group $\G$, it is also natural to investigate how an isometric action by $\G$ on a manifold $M \in  \mathcal{P}$ restricts $M$. For example, in which range of dimensions can one get a complete classification, or an exhaustive list of possibilities, potentially containing new examples. Also

\begin{problem}
Is a simply connected $\G$ manifold $M^n \in \mathcal{P}$ of minimal dimension a CROSS with a linear action?
\end{problem}

In several of the above, as well as in subsequent results, the classification of \emph{fixed point homogeneous} manifolds of positive curvature by Grove-Searle plays an important role. Here a manifold is called fixed point homogeneous if it supports an action where the group acts transitively on the normal sphere to a fixed point component, i.e., cofix = 0 with the terminology provided above. A similar useful classification is also known in fixed point cohomogeneity one by Grove-Kim:

\begin{thm}[Almost Minimal Cofix]
Any  simply connected $M \in \mathcal{P}$ with ${\ensuremath{\operatorname{cofix}}} (M) \le 1$  is diffeomorphic to a CROSS.
\end{thm}

There is also a complete classification up to equivariant diffeomorphism in the non simply connected case. The conclusion fails in cofix 2.

\bigskip

For the degree of symmetry one has the following satisfactory result due to Wilking:

\begin{thm}[Degree of Symmetry]
Any simply connected $M^n \in \mathcal{P}$ with $\symdeg (M) \ge 2n - 6$ is tangentially homotopy equivalent to a CROSS, or isometric to a homogeneous space in $\mathcal{P}$. 
\end{thm}

One of the striking things about the classification of homogeneous manifolds of positive curvature is that apart from the CROSS'es, they occur only in finitely many dimensions. This phenomenon turns out to hold for any cohomogeneity by the following result of Wilking:

\begin{thm}[Cohomogeneity Finitenes]
Any simply connected $M^n \in \mathcal{P}$ with $\cohom(M) = k \ge 1$ and $n \ge 18(1+ k)^2$ is tangentially homotopy equivalent to a CROSS.
\end{thm}

One of the remarkable original constructions, which lead to this result, was the existence of a \emph{chain} of positively curved $G_i$ manifolds  $M_i$ (with fixed codimension) and with isometric orbit spaces $M_i/G_i$ once a sufficiently large dimensional positively curved cohomogeneity  $k$ manifold is given. A deep analysis of the limit object $M_{\infty}$ showed that it was one of $\Sph^{\infty}, \CP^{\infty}$ or $\HP^{\infty}$.

\smallskip

While the above result illustrates the difficulty in finding new examples of manifolds in $ \mathcal{P}$, it also gives hope that it might be possible for low cohomogeneity. 

\smallskip

In the case of cohomogeneity one, the bound in the above theorem is 72. The optimal number is actually 14. In fact, from the classification in even dimensions due to Verdiani and the work of Grove-Wilking-Ziller, one has:

\begin{thm}[Cohomogeneity One]
A simply connected cohomogeneity one $n$-manifold $M \in \mathcal{P}$ is diffeomorphic to a CROSS when $n \ne 7, 13$.

For $n = 13$, $M$ is diffeomorphic to a cohomogeneity one Bazaikin space or $\Sph^{13}$.

For $n=7$, $M$ is diffeomorphic to a cohomogeneity one Eschenburg space, the Berger space $B^7$, or to one of the Konishi-Hitchin manifolds $P_k$, $Q_k$, $k \ge 2$, or to an exceptional manifold $N^7$.
\end{thm}

In all cases, the possible actions are known as well (for the CROSS'es they are all ``linear"). One of the Bazikin manifolds in the theorem is the Berger space $B^{13}$, and the theorem gives a classification in this dimension. One of the Eshenburg spaces and $Q_1$ is the Aloff-Wallach-Wilking normal homogeneous manifold $A_{11} = W^7 = \SU(3) \SO(3)/\U(2)$, and $P_1 = \Sph^7$. The theorem provides a classification of positively curved cohomogeneity one manifolds modulo the now obvious

\begin{problem}
Determine whether $P_k$, $Q_k$, $k \ge 2$ and $N$ support (invariant) metrics of positive curvature.
\end{problem}

All the 7-manifolds, $M$ in the above theorem support almost effective actions by $\S^3 \times \S^3$ and metrics of nonnegative curvature (Grove-Ziller). Moreover, each of the $\S^3$ factors act almost freely on $M$. In particular, there is an \emph{orbifold} bundle $M \to M/\S^3 = B$, where $B$ supports an almost effective action by the other $\S^3$ factor. In fact, $B = \Sph^4$ or $\CP^2$ and the action is the well known action with singular orbits of codimension 2. For $B^7$ and $N^7$ the metric on the base is an orbifold metric, singular along both singular orbits. For $P_k, Q_k, k \ge 2$ the metric on the base is an orbifold metric which is singular along one of the singular orbits.
In the latter cases Hitchin has constructed a self dual einstein metric on each base, and the manifolds $P_k$ and $Q_k$ are the two-fold covers of the so-called Konishi bundles of self dual two forms of these. By work of Dearricott it is well known that if the Hitchin metrics had positive curvature then the natural 3-Sasakian metric on the Konishi bundle would also have positive curvature once the fiber is shrunk sufficiently. Alas, the Hitchin metrics have curvatures of both signs when $k \ge 2$.

The most interesting among the candidates are the $P_k$'s since they are all 2-connected. In particular, it would follow from the $\pi_2$ finiteness theorem (see Shape section), that if they have positive curvature, their pinching must necessarily approach 0 as $k$ goes to infinity. Although these candidates have been known for several years, some fundamental new ideas are still needed to see if this family as a whole supports positive curvature. However,

\begin{thm}[Cohomogeneity One Example]
The manifold $P_2$ has an invariant metric with positive curvature.
\end{thm}

This indeed is a new example, since $P_2$ is $2$-connected with $\pi_3(P_2) = \Z_2$, and the only known two connected 7-manifolds with positive curvature are $\Sph^7$ and $B^7$ and $\pi_3(B^7) = \Z_{10}$. It is interesting to note that in fact $P_2$ is homeomorphic to the unit tangent bundle $T_1\Sph^4$ of the 4-sphere, but it is not known if it is diffeomorphic to it.

As indicated earlier, any construction of a metric with positive curvature on any of the new candidates must be unlike all previous ones. For \emph{connection metrics} a necessary and sufficient condition for positive curvature has been derived by 
Chaves-Derdzi\'{n}ski-Rigas. In a manuscript \cite{De}, Dearricott offers a proof that this condition can be satisfied by means of a suitable conformal change of the  Hitchin metric. A very different construction of an invariant piecewise polynomial $C^2$ metric with positive curvature has been given by Grove-Verdiani-Ziller in \cite{GVZ}. In fact, that metric is shown to have \emph{strongly positive curvature}, i.e., there is an (explicitly constructed) 2-form $\eta$ such that the modified curvature operator $\hat R + \hat \eta$ (necesseraly having the same ``sectional curvatures") is positive definite. Checking this amounts to checking that specific polynomials are positive on a given interval, the orbit space of the group action.

\bigskip

Another natural attempt for the construction of new examples of positive curvature, is to see if one or more among the large known class of nonnegatively curved manifolds can be \emph{deformed} to have positive curvature. Although no obstructions are known yet for simply connected manifolds, this is exceedingly difficult. By the orientation theorem of Synge, $\RP^n \times \RP^m$ provide simple examples where this is impossible. Here $\RP^3 \times \RP^2$ is particularly striking, since Wilking has constructed  a metric with nonnegative curvature, and positive curvature on an open dense set of points. The same property is known for the so-called Gromoll-Meyer sphere  $\Sigma^7= \Sp(2)//\Sp(1)$ by work of Wilhelm. In a recent manuscript \cite{PW},  Petersen and Wilhelm  offer a proof that for this example a deformation is possible:
\smallskip

\begin{thm}[Exotic Sphere]
$\Sigma^7= \Sp(2)//\Sp(1)$ admits a metric with positive curvature.
\end{thm}

\smallskip

It is remarkable that during the deformation, metrics with some negative curvatures will arise and the whole proof is both long and involved. The existence of an exotic sphere with positive curvature is of course of pivotal importance for all differentiable sphere theorems. It should also be noted that by work of Weiss and Grove-Wilhelm respectively, $\Sigma^7$ does not support a 1/4 pinched (cf. Shape Section) metric (nor one with \emph{radius} larger than 1/2 maximal) , respectively a positively curved metric with 4 points at distance exceeding 1/2 of the maximal possible diameter (cf. Size Section).

\bigskip

We like to mention that at the time of writing, the work on all of the above examples has not yet been through a complete examination by other experts.

\bigskip

Note that from the classification theorem above we know in particular that any cohomogeneity one homotopy $n$-sphere $\Sigma \in \mathcal{P}$ is standard. The same conclusion holds in the class $\mathcal{P}_0$ by work of Grove-Verdiani-Wilking-Ziller, but not in $\mathcal{P}_0^{-}$. In fact all cohomogeneity one manifolds belong to $\mathcal{P}_0^{-}$ by Schawachh\"{o}fer-Tuschmann, but the Kervaire spheres have cohomogeneity one and some of them are exotic. Since the proposed metric of positive curvature in the Petersen-Wilhelm example  has cohomogeneity 4, and a positively curved homotopy 4-sphere with a circle action is standard by the strong rank theorem, the following is natural:

\begin{problem}
Is a positively curved cohomogeneity $k$ homotopy sphere standard if $k \le 3$?
\end{problem}

 It is also natural to wonder about whether or not actions are linear. As we have seen they are linear in cohomogeneity one for the class $\mathcal{P}$, and almost for the class $\mathcal{P}_0$. The exception here follows from work Grove-Ziller, since there is a non-linear ``Kervaire" action, observed by Calabi on $\Sph^5$ where both singular orbits have codimension 2, and hence the manifold has an invariant metric with non-negative curvature.

\section{Size}

There are several natural metric invariants measuring \emph{size}, such as, e.g., \emph{diameter, radius} and \emph{volume}. This allows for precise ways of saying that positively curved manifolds are ``small". The first of these, obtained by a Synge type argument,  due to Bonnet (for sectional curvature) and Myers (for Ricci curvature) provides a diameter bound in terms of a positive lower curvature bound:

\begin{thm}[Diameter Bound]
A manifold with curvature bounded below by $\delta > 0$ has diameter at most $\pi/\sqrt{\delta}$, the diameter of the sphere with constant curvature $\delta$.
\end{thm}

In particular, a positively curved manifold has finite fundamental group, and in even dimensions it is either trivial or $\Z_2$ by Synge's theorem. No further general restrictions on the structure of the fundamental group of a positively curved manifold are known (cf. Shape section). In particular, the first \emph{counter example} to the so-called \emph{Chern conjecture} proposing that an abelian fundamental group would be cyclic, was found by Shankar. It is easy to see that any finite group is a subgroup of some $\SU(n)$, and hence the fundamental group of a manifold in $\mathcal{P}_0$. In contrast, nothing is known about the natural:

\begin{problem}
Is any finite group the fundamental group of a positively curved manifold?
\end{problem}

\no In any given dimension, however, there are restrictions on the \emph{size of the topology} provided by the following finiteness result  due to Gromov:

\begin{thm}[Betti Number Theorem]
For any $n \in \N$ there is a $C(n)$ such that the fundamental group and homology groups of any $M^n \in  \mathcal{P}_{0^{-}}$ is generated by at most $C(n)$ elements.
\end{thm}

The main geometric tools used in the proof of this result, is the Toponogov Comparison Theorem, critical point theory for distance functions and the Bishop-Gromov \emph{Relative Volume Comparison Theorem}. 

\begin{problem}
Determine the optimal value for $C(n)$.
\end{problem}

An extremely strong size and structure restriction is contained in the following proposed extension of a conjecture attributed to Bott:

\begin{conj}[Ellipticity]
The Betti numbers of the loop space (any field of coefficients) of any simply connected $M^n \in  \mathcal{P}_{0^{-}}$ grow at most polynomially.
\end{conj}

A confirmation would have profound consequences. When restricted to the field of rational numbers a lot of restrictions are known for this class of so-called \emph{rationally elliptic manifolds}. In particular, it would yield the optimal bound $C(n) = 2^n$ for the sum of Betti numbers, and also have impact on the classical Hopf conjecture about the Euler characteristic. The extension to the class $\mathcal{P}_{0^{-}}$ could be significant in a potential path of a proof of the ellipticity conjecture.

\bigskip

The Bonnet-Myers theorem begs the question of what happens in the extreme case. In fact, having \emph{maximal diameter characterizes the round sphere} as was shown by Toponogov:

\begin{thm}[Maximal Diameter Theorem]
A manifold with curvature bounded below by $\delta > 0$ and diameter $\pi/\sqrt{\delta}$ is isometric to the sphere with constant curvature $\delta$.
\end{thm}

The conclusion also holds when sectional curvature is replaced by Ricci curvature.

\bigskip

For positively curved manifolds, say with lower curvature bound normalized to $1$, the \emph{diameter function} therefore provides a natural filtration in terms of super level sets, to be thought of as ``diameter pinching". In particular, it is natural to ask whether a sufficiently large manifold with positive curvature is a sphere. Rather than relying on classical Morse theory for the loop space, it turned out to be advantageous to work directly on the manifold and develop a ``Morse type" theory for (non-smooth) distance functions. In conjunction with Toponogov's comparison theorem, this indeed lead to the following results of Grove-Shiohama, and Gromoll-Grove and Wilking

\begin{thm}[Diameter Sphere Theorem]
A manifold with curvature bounded below by $\delta$ and diameter  $ > \pi/2\sqrt{\delta}$ is homeomorphic to a sphere.
\end{thm}

\begin{thm}[Diameter Rigidity Theorem]
A manifold with curvature bounded below by $\delta$ and diameter $\pi/2\sqrt{\delta}$ is either a topological sphere, or its universal cover is isometric to a CROSS.
\end{thm}

In the latter case, there is also an isometric classification when the fundamental group is non-trivial: The $\Z_2$ quotient of complex odd dimensional projective spaces, and  the space forms where the fundamental group acts irreducibly. 

\smallskip

These ``diameter theorems" raise two natural questions: 

\begin{problem}
Is a manifold in $\mathcal{P}$ with almost maximal diameter diffeomorphic to the standard sphere?
\end{problem}

\no or even

\begin{problem}
Is a manifold in $\mathcal{P}$ with almost 1/2 maximal diameter diffeomorphic to one of the model spaces?
\end{problem}

Problems of this type are naturally analyzed using the fact that the class of $n$-manifolds with given lower bound on Ricci curvature and upper bound on diameter is \emph{precompact} relative to the so-called \emph{Gromov-Hausdorff metric}. This is a fairly simple consequence of the relative Bishop-Gromov volume comparison theorem. When a lower (sectional) curvature bound is present, any limit object is a so-called \emph{Alexandrov space}, i.e., a finite dimensional \emph{length/inner metric} space with a lower curvature bound expressed by \emph{distance comparison} equivalent for manifolds to the Toponogov triangle comparison theorem. The dimension of such an Alexandrov space is at most $n$. One refers to \emph{collapse} in the general case of dimension less than $n$, and \emph{non-collapse}, otherwise. Non-collapse is equivalent to having a lower bound on volume. Relatively little is known in general about collapse. The work so far has culminated in the recent manuscript \cite{KPT} by Kapovitch-Petrunin-Tuschman  on the class $\mathcal{P}_{0^{-}}$. In low dimensions, much more is known due to work of Shioya-Yamaguchi which has significant impact on the solution of Thurston's \emph{geometrization conjectuture} due to Perelman. In the non-collapsing case, however, one has the following fundamental result due to Perelmann, a proof of which has recently been published in \cite{Ka}:

\begin{thm}[Topological Stability]
All $n$-dimensional Alexandrov spaces with curvature $\ge k$ in a Gromov-Hausdorff neighborhood of an $n$-dimensonal Alexandrov space $X$ with curvature $\ge k$ are homeomorphic to $X$.
\end{thm}

This immediately yields topological finiteness in all dimensions (including 3 where only homotopy type was know) for the class of $n$-manifolds with lower bounds on curvature and volume and upper bound on diameter, originally due to Grove-Petersen-Wu. In the case of positive curvature (where also Hamilton's theorem can be invoked), one gets in particular:

\begin{thm}[Topological Finiteness]
For any positive $\delta, v$ and integer $n\ne 4$, there are at most finitely many  diffeomorphism types of $n$-manifolds, with curvature and volume bounded below by $\delta$ and $v$ respectively. The same holds for homeomorphism types in dimension 4.
\end{thm}

Note, that in the diameter problems above no volume bound is present, so collapse will typically happen. Of course limiting objects are Alexandrov spaces with positive curvature and maximal, respectively half maximal diameter. In the case of maximal diameter, the Alexandrov space is not a round sphere as in Toponogov's theorem above, but ``only" a so-called \emph{spherical suspension} of an Alexandrov space with the same positive lower curvature bound.

The spherical suspension of a manifold in $\mathcal{P}$ with more than 1/2 maximal diameter was closely analyzed by Grove-Wilhelm, and used to derived the following result where collapse is avoided:

\begin{thm}
For any $n$-manifold $M \in \mathcal{P}$ with more than $1/2$ maximal diameter, there is an $n$-dimensional Alexandrov space $X$ which is the Gromov-Hausdorff limit of positively curved metrics on $M$ as well as on the standard $n$-sphere. 
\end{thm}

In other words, a \emph{differentiable stability} version of Perleman's Theorem above would imply that such an $M$ is diffeomorphic to the standard sphere. A curious observation by Wilking based on the suspension construction as well, also shows that if there is an exotic sphere $M \in \mathcal{P}$ with more than $1/2$ maximal diameter, then $M$ also carries a metric with almost maximal diameter. In other words, the answer to problem 11 is yes if and only if any positively curved manifold with diameter larger than 1/2 maximal is diffeomorphic to a sphere.

The case where $M \in \mathcal{P}$ has $q+1 > 2$ points with individual distances larger than $1/2$ maximal diameter was also analyzed  via spherical suspensions. It was shown, that $M = \Disc^q \times \Sph^{n-q} \cup \Sph^{q-1} \times \Disc^{n-q+1}$. In particular, $M$ is diffeomorphic to the standard sphere if $q = n-3$ invoking a result of Hatcher.

\bigskip

All size invariants are difficult to compute, and so filtrations by size do not provide any guide towards the construction of new examples.

\section{Shape}

A third \emph{characterization of the round sphere} is being simply connected and having \emph{constant positive curvature}. It is only natural to expect that (simply connected) Riemannian manifolds with almost the same shape are spheres as well. This became known as the \emph{pinching problem} raised by Hopf. Here, the \emph{pinching}

\begin{center}
$\delta_M := \frac{\min sec}{\max sec} \le 1$
\end{center}

\no measures the proximity to constant curvature. The first major result in this direction was due to Rauch who in the process derived the so-called \emph{Rauch comparison theorem} for Jacobi fields. The optimal result in terms of pinching, the so-called $1/4$- pinching \emph{sphere theorem},  was proved by Klingenberg, with \emph{rigidity} due to Berger:

\begin{thm}[Topological Sphere/Rigidity Theorem]
A $\delta$ pinched simply connected manifold with $\delta \ge 1/4$ is either homeomorphic to the sphere or isometric to a CROSS.
\end{thm}

  The crucial geometric tools used to prove this was the Toponogov comparison theorem and the so-called \emph{long homotopy lemma} by Klingenberg. The question of diffeomorphism was first succesfully treated independently by Calabi, Gromoll and Shikata each of whom established a diffeomorphism sphere theorem with dimension dependent pinching approaching 1. Each proof also provided the germ for subsequent developments: The \emph{Gromoll filtration} of the group of homotopy spheres, the idea and notion of \emph{distance between manifolds/differentiable structures} in the work of Shikata, and the analytic approach used by Calabi.  In particular, some of the roots of the far reaching work of Colding and Cheeger on manifolds with a lower Ricci curvature bound have similarities with the work of Calabi. - Subsequently, two different methods to achieve a dimension independent diffeomorphism sphere theorem were obtained by  Shiohama-Sugimoto, and by Ruh. Moreover, the method of Ruh in conjunction with the general non-linear \emph{center of mass} technique developed by Grove-Karcher led to an equivariant sphere theorem, and in particular to a pinching theorem for all \emph{space forms}.
 
 Just recently, following a breakthrough by B\"{o}hm and Wilking \cite{BW}, the \emph{Ricci flow} was used by Brendle and Schoen  \cite{BS1,BS2} to replace homeomorphism by diffeomorhism in the above classical result, and extend it to pointwise pinching and equivariance. In particular:
 
 \begin{thm}[Differentiable Sphere/Rigidty Theorem]
A pointwise 1/4-pinched manifold is either diffeomorphic to a space form or locally isometric to a CROSS
\end{thm}

A key issue here is which curvature conditions are preserved by the Ricci flow, a very difficult problem. Although $1/4$ pinching is not preserved, it turns out that it is in a family which is preserved and which moreover satisfies the general condition set forth by 
B\"{o}hm and Wilking in the situation of strict 1/4 pinching.

\smallskip

Already some time ago, Berger using Gromov-Hausdorff limit arguments obtained a below 1/4 pinching theorem in even dimensions (for spheres only deriving homeomorphism of course), and later using non convergence methods, Abresh-Meyer were able to do the same in odd dimensions getting explicit constants. Invoking the Ricci flow, Petersen-Tao \cite{PT} recently derived the satisfactory

 \begin{thm}[Below 1/4 pinching]
There is an $\epsilon(n) > 0$ such that any simply connected  $1/4 - \epsilon(n)$ pinched $n$-manifold is diffeomorphic to a CROSS.
\end{thm}

The first $\delta < 1/4$ where one knows for sure that there is a different $\delta$-pinched manifold, is for $\delta = 1/37$. In fact 
P\"utmann showed that each of the three normal homogeneous manifolds, $B^7 = \SO(5)/\SO(3), W^7 = \SU(3)\SO(3)/\U(2)$, and $B^{13} = \SU(5)/\Sp(2)\S^1$ with positive curvature admits a homogeneous 1/37-pinched metric. It is natural to wonder if indeed $1/37$ is the \emph{optimal pinching} among all positively curved metrics on these examples, and how large $\epsilon(n)$ necessarily has to be. 
\bigskip

An investigation of the class of $\delta$-pinched manifolds equipped with the Gromov-Hausdorff metric does provide some crude insight in terms of  \emph{finiteness results} or even obstructions modulo finiteness.

\smallskip

The general \emph{Cheeger Finiteness Theorem} asserts that the class of closed $n$-manifolds with bounded curvature and diameter, and with a lower bound for volume contains at most finitely many diffeomorphism types. Since a positive lower curvature bound by the Bonnet-Myers Theorem yields a bound for the diameter, and an upper curvature bound yields a lower bound on the injectivety radius for positively curved manifolds in even dimensions by results of Synge and Klingenberg, one has the following 

 \begin{thm}
For any $n \in \N$, and $0 < \delta \le 1$, there are only finitely many diffeomorphism types of $\delta$-pinched $2n$-manifolds.
\end{thm}

So far, at most 4 different simply connected examples of positively curved $2n$-manifolds are known for any $n$.
\smallskip

In contrast, one easily sees that there exists a sequence of different Aloff-Wallach spaces $A_{pq}$
whose pinching converges to that
of $W^7 = A_{11}$ and hence to 1/37. In particular, for any $\delta < 1/37$ there are infinitely many $\delta$-pinched manifolds. Nevertheless, with additional topological restrictions, the following is a corollary of $\pi_2$-\emph{Finiteness Theorem} due to Petrunin-Tushmann and to Fang-Rong:

 \begin{thm}
For any $n \in \N$, and $0 < \delta \le 1$, there are only finitely many diffeomorphism types of $\delta$-pinched 2-connected $n$-manifolds.
\end{thm}

So far, however, no infinite family of positive curvature is known, where by necessity the pinching constants will approach $0$. An infinite family in even dimensions, or an infinite 2-connected family in odd dimensions would do according to these finiteness results. The $P_k$ family from the Symmetry section is a natural candidate for this. On the other hand, Fang-Rong have conjectured that the above theorem should hold for positively curved manifolds without a pinching assumption.

Utilizing the Cheeger-Fukaya-Gromov theory of collapse with bounded curvature, Rong has obtained he following result about the structure of the fundamental group

 \begin{thm}[$\pi_1$ Structure]
For any $n \in \N$, and $0 < \delta \le 1$, there is a constant $w(n,\delta)$ so that $\pi_1(M)$ contains a cyclic subgroup of index at most $w(n,\delta)$ for any $\delta$-pinched $n$-manifold $M$.
\end{thm}

It is a conjecture of Rong that $w(n,\delta)$ can be chosen independently of $\delta$.

\bigskip

Although the pinching function provides a natural filtration among positively curved manifolds by superlevel sets, it is typically very difficult to calculate for a particular example and clearly not a guide towards finding new examples.

 \providecommand{\bysame}{\leavevmode\hbox
to3em{\hrulefill}\thinspace}

\end{document}